\documentclass[anonymous]{llncs}
\usepackage{graphicx}
\usepackage[utf8]{inputenc}
\usepackage{amsmath}
\usepackage{amsfonts}       
\usepackage{booktabs}    
\usepackage{xcolor}      
\usepackage{hyperref}
\usepackage{todonotes}
\usepackage[]{lineno}
\usepackage{orcidlink}

\newcommand*\samethanks[1][\value{footnote}]{\footnotemark[#1]}

\title{Reconfigurations of Plane Caterpillars and Paths}

\author{
Todor Anti\'c \inst{1}\thanks{Supported by project 23-04949X of the Czech Science Foundation
(GA\v{C}R) and by scholarships provided by the Visegrad Fund.}\orcidlink{0009-0008-6521-7987}
\and
Guillermo Gamboa Quintero \inst{2}\thanks{Supported by the GA\v{C}R grant 22-17398S.}\orcidlink{0000-0002-8968-6269}
\and
Jelena Gli\v{s}i\'c \inst{1}\samethanks[1]\orcidlink{0009-0002-0792-3070}
}
\institute{Department of Applied Mathematics, 
Charles University, Prague, Czechia \email{\{todor,glisic\}@kam.mff.cuni.cz} \\
\and 
Computer Science Institute of Charles University, Prague, Czechia \email{gamboa@iuuk.mff.cuni.cz}
}

\authorrunning{T. Anti\'c, G. Gamboa Quintero and J. Gli\v{s}i\'{c}}

\newcommand{\starat}{\operatorname{star}}

\graphicspath{{figures/}}

\definecolor{dark blue}{rgb}{0.121,0.47,0.705}
\let\emph\relax\DeclareTextFontCommand{\emph}{\color{blue}\em}

\begin{document}

\maketitle

\begin{abstract}
    Let $S$ be a point set in the plane, $\mathcal{P}(S)$ and $\mathcal{C}(S)$ sets of all plane spanning paths and caterpillars on $S$. We study reconfiguration operations on $\mathcal{P}(S)$ and $\mathcal{C}(S)$. In particular, we prove that all of the commonly studied reconfigurations on plane spanning trees still yield connected reconfiguration graphs for caterpillars when $S$ is in convex position. If $S$ is in general position, we show that the rotation, compatible flip and flip graphs of $\mathcal{C}(S)$ are connected while the slide graph is disconnected. 
    For paths, we prove the existence of a connected component of size at least $2^{n-1}$ and that no component of size at most $7$ can exist in the flip graph on $\mathcal{P}(S)$.
    \keywords{Reconfiguration graph, Geometric graph, Caterpillar, Path}
\end{abstract}

\section{Introduction}
 Given a set of structures $C$, and a reconfiguration operation that transforms one object in $C$ to another, the \emph{reconfiguration graph} is a graph with vertex set $C$ in which two vertices form an edge if one can be transformed into the other using a single reconfiguration operation. Often, in computer science, objects are solutions to a problem and reconfigurations are local changes that transform one solution into another. Then, to understand the solution space of a problem, it is important to study both the structural properties of the reconfiguration graph (connectivity, hamiltonicity, etc.) and algorithmic questions (how to find the shortest reconfiguration sequence). For an introduction to the topic of reconfiguration, see \cite{Nishimura2018survey}. We  focus on reconfigurations in the following setting.
 
 Given a point set $S$ in the plane, a \emph{plane spanning tree} on $S$ is a spanning tree of $S$ whose edges are straight line segments that do not cross. Let $\mathcal{T}(S)$ be the set of all plane spanning trees on $S$. We define the following five reconfigurations on $\mathcal{T}(S)$. For the following, we are given plane spanning trees $T_1 = (S,E_1), T_2 = (S,E_2) \in \mathcal{T}(S)$, then we say that: 

 \begin{enumerate}
     \item $T_1$ and $T_2$ are connected by a \emph{flip} if $E_2 = E_1 \setminus \{e\} \cup \{f\}$ for some edges $e,f$.
     \item $T_1$ and $T_2$ are connected by a \emph{compatible flip} if $E_2 = E_1 \setminus \{e\} \cup \{f\}$ for some edges $e,f$ which do not cross.
     \item $T_1$ and $T_2$ are connected by a \emph{rotation} if $E_2 = E_1 \setminus \{e\} \cup \{f\}$ for some edges $e ,f$ which share an endpoint.
     \item $T_1$ and $T_2$ are connected by an \emph{empty triangle rotation} if $E_2 = E_1 \setminus \{e\} \cup \{f\}$ for some edges $e,f$ which share an endpoint and the triangle spanned by their endpoints is empty.
     \item  $T_1$ and $T_2$ are connected by a \emph{slide} if $E_2 = E_1 \setminus \{e\} \cup \{f\}$ for some edges $e,f$ which share an endpoint and the triangle spanned by their endpoints is empty and if $e= ab$ and $f=ac$ then $bc \in E_1 \cap E_2$. 
 \end{enumerate}

\begin{figure}
     \centering
     \includegraphics[width=1\linewidth]{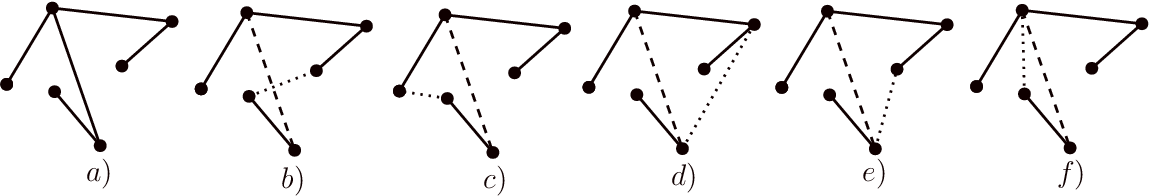}
     \caption{a) A plane spanning tree. Reconfiguration of the tree by changing the dashed line to the dotted line is: b) a flip, c) a compatible flip, d) a rotation, e) an empty triangle rotation and f) a slide.}
     \label{Fig: rekonfiguracije}
 \end{figure}

 For a visualization of all of the reconfiguration operations described above, see Figure \ref{Fig: rekonfiguracije}. From the description of these operations, one can notice that there exists a linear hierarchy. Every slide is an empty triangle rotation, every empty triangle rotation is a rotation, and so on. This hierarchy is useful when studying the structural properties of the corresponding reconfiguration graphs. For example, if the \emph{slide graph of plane spanning trees} is connected then so are all of the other reconfiguration graphs. 

 The reconfiguration graphs associated with the operations described above have been a topic of interest for a long time with many results appearing through the years. These results have concerned connectivity \cite{slideconn,TreeEnumeration,NICHOLS2020111929}, lower and upper bounds on the diameter \cite{TreeGraphDiameter,bousquet_et_al:LIPIcs.SoCG.2024.22,TreeConvexDiameter}, etc. However, there has been little research on induced subgraphs of these reconfiguration graphs. The only such subgraph that has been explored is the subgraph induced by plane spanning paths. And even then, some of the main questions have been for a long time. We aim to expand on the study of induced subgraphs of reconfiguration graphs of plane spanning trees by exploring the previously unexplored problem of reconfigurations of plane spanning caterpillars, and by expanding on the topic of reconfigurations of plane spanning paths.

\subsection{Our contribution}

A \emph{caterpillar} is a tree in which all non-leaf vertices form a path. Possibly, this path is a single vertex or empty if the caterpillar is a single edge. We call this path the \emph{spine} of the caterpillar. For a point set $S$ in the plane, a \emph{plane spanning caterpillar} of $S$ is a plane spanning tree of $S$ which is a caterpillar. We call the two endpoints of the spine the \emph{head} and the \emph{tail} (we may choose which one is which) and leaves connected to these vertices \emph{head-leaves} and \emph{tail-leaves}. For a set $S$, we will denote by $\mathcal{C}(S)$ the set of all plane spanning caterpillars on $S$. We will denote the reconfiguration graphs on $\mathcal{C}(S)$ by $G_{\mathcal{C}}^{\text{flip}}(S),G_{\mathcal{C}}^{\text{comp-flip}}(S),G_{\mathcal{C}}^{\text{rot}}(S),G_{\mathcal{C}}^{\text{emp-rot}}(S),G_{\mathcal{C}}^{\text{slide}}(S).$

 First, we focus on the case where $S$ is in convex position. We show that the \emph{slide graph }$G_{\mathcal{C}}^{\text{slide}}$ is connected in this case.

\begin{theorem}\label{Thm: Convexpillars}
    Let $S$ be a set of $n\ge3$ points in convex position in the plane. Then, the graph $G_{\mathcal{C}}^{\text{slide}}(S)$ is connected with diameter at most $3n-8$.
\end{theorem}

Of course, this immediately implies that all of the other reconfiguration graphs are connected for a point set $S$ in convex position.  Then, we consider the case where $S$ is a point set in general position. Here the situation is very different. Mainly, for each $n\ge 9$, it is possible to find sets $S$ of $n$ points such that $G_{\mathcal{C}}^{\text{slide}}(S)$ has isolated vertices. On the other hand, we can prove connectivity for the \emph{rotation graph} $G_{\mathcal{C}}^{\text{rot}}(S)$. 

\begin{proposition}\label{prop: rotations}
Let $S$ be a set of points in general position in the plane. Then, the graph $G_{\mathcal{C}}^{\text{rot}}(S)$ is connected.   
\end{proposition}

The above proposition again implies that the \emph{flip graph}  $G_{\mathcal{C}}^{\text{flip}}(S)$ and \emph{compatible flip graph} $G_{\mathcal{C}}^{\text{comp-flip}}(S)$ of plane spanning caterpillars are connected for each $S$ in general position. However, our results do not imply connectivity of the \emph{empty rotation graph} $G_{\mathcal{C}}^{\text{emp-rot}}(S)$, so the following question remains open. 

\begin{question}
    Is $G_{\mathcal{C}}^{\text{emp-rot}}(S)$ connected for $S$ in general position?
\end{question}

Given the disconnectedness of $G_{\mathcal{C}}^{\text{slide}}(S)$, it becomes more interesting to find large connected components in this graph. To do this, we consider special subclasses of caterpillars and show that these are connected in $G_{\mathcal{C}}^{\text{slide}}(S)$. We write $\starat(x)$ for a spanning star with center $x$. 

\begin{lemma}\label{Thm: starsanddoublestars}
    Let $S$ be a set of $n$ points in general position in the plane and $u,v\in S$ be adjacent in the convex hull of $S$. Then, $\starat(u)$ and $\starat(v)$ are both connected to any double star with centers $u$ and $v$ in $G_{\mathcal{C}}(S)$. 
\end{lemma}

For a caterpillar $C\in \mathcal{C}(S)$, and consecutive spine vertices $v_i,\dots, v_j$ of $C$, we write $S_{i,j}$ for the point set consisting of the spine vertices and all of the leaves attached to them.
We call $C\in \mathcal{C}(S)$ with spine $v_1,v_2,\dots, v_k$ a \emph{well-separated} caterpillar if for each $i\ge 1$, the convex hull of $S_{1,i}$ is disjoint from the rest of $S$. As a consequence of Lemma \ref{Thm: starsanddoublestars}, we get that all caterpillars in this relatively general class are mutually connected in $G_{\mathcal{C}}(S)$. 

\begin{theorem}\label{Thm: wellseparated}
     Any two well-separated caterpillars are connected by a sequence of slides in $G_{\mathcal{C}}(S)$.
\end{theorem}

Then, we solve the connectivity of $G_{\mathcal{C}}^{\text{slide}}(S)$ for $S$ in general position. 

\begin{theorem}\label{thm: sliidesconnnected}
    Let $n$ be a natural number. Then $G_{\mathcal{C}}^{\text{slide}}(S)$ is connected for every set $S$ of $n$ points in the plane if $n\le 7$. If $n\ge 8$, there exists a set $S$ of $n$ points such that $G_{\mathcal{C}}^{\text{slide}}(S)$ has isolated vertices. 
\end{theorem}

Lastly, we shift our focus to the study of connected components of the reconfiguration graph of plane spanning paths.  Given a set of points in general position $S$, we will refer to the corresponding \emph{flip graph of plane spanning paths} as $G_{\mathcal{P}}(S)$. Currently, the main open problem related to such flips is deciding if $G_{\mathcal{P}}(S)$ is connected. In this direction, we prove the following result. 

\begin{theorem}\label{Thm: largecomponent}
    Let $S$ be a set of $n$ points in general position. Then $G_{\mathcal{P}}(S)$ contains a connected component of size $\Omega(2^{n-1})$.
\end{theorem}

In particular, this component consists of the paths which we call \emph{generalized peeling paths}. We introduce this subclass in Section \ref{Sec: flips}.  

We note that Theorem \ref{Thm: largecomponent} was independently discovered by Kleist, Kramer and Rieck \cite{personalcom,Lindapaper}. We still include it here because we use the number of generalized peeling paths to prove that there are at least $\frac{1}{2}(3^{n}-1)$ well-separated caterpillars which implies that $G_{\mathcal{C}}^{\text{slide}}(S)$ has a component of at least this size. 

Finally, we investigate the minimal size of components in $G_{\mathcal{P}}(S)$. In this direction, we prove the following result. 

\begin{theorem} \label{thm: smallcomponents}
Let $S$ be a point set of $n\ge 5$ points in general position. Then, $G_{\mathcal{P}}(S)$ contains no connected component of size at most $7$.     
\end{theorem}

\subsection{Previous work}

\subsection*{Reconfiguration graphs of plane spanning trees}
It is known that the reconfiguration graph of plane spanning trees is connected even in the most restrictive case when the reconfigurations are slides \cite{slideconn}. Consequently, so are the reconfiguration graphs for all other types of reconfigurations we have defined. In the case of slides, a tight $\Theta(n^2)$ bound on the diameter is shown in \cite{quadratic}. For the empty triangle rotations, an upper bound of $O(n\log n)$ on the diameter was shown in \cite{NICHOLS2020111929}, while for the remaining reconfigurations a linear upper bound is known \cite{TreeEnumeration}. In \cite{TreeGraphDiameter}, an upper bound of $2n-3$ is shown for flip graphs. With the exception of slide graphs, the best known lower bound on the diameter of the reconfiguration graphs is $1.5n-5$, as per \cite{slidelb}.

\subsection*{Reconfiguration graphs of plane spanning paths}

The reconfiguration graph of plane spanning paths has been thoroughly studied for convex point sets. It is known that the flip graph is connected \cite{pathflipsconv} and that for a convex point set of size $n$ the flip graph has diameter $2n-5$ for $n \in \{ 3, 4 \}$ and $2n-6$ for $n \geq 5$ \cite{pathflipsconvdiam}. Moreover, the flip graph is Hamiltonian \cite{PathGraphIsHamiltonian} and it has chromatic number $n$ \cite{chromnumberconvex}. Also, it is known that the flip graph of a convex point set of size $n$ where the paths considered have a fixed start vertex has diameter $2n-5$ and radius $n-2$ for $n \geq 3$ \cite{Lindapaper}. Finally, the flip graph of a point sets with at most two convex layers is connected \cite{Lindapaper}. For point sets in general position connectivity was first conjectured by Akl et al. \cite{pathflipsconv}. 

\begin{conjecture}[Akl. et. al. \cite{pathflipsconv}]\label{conj: main}
Let $S$ be a point set in general position. Then the flip graph of plane spanning paths on $S$ is connected.
\end{conjecture}

Despite this conjecture being around 17 years old, there has been relatively little progress towards solving it. Connectivity has been shown for point sets of size $n \leq 8$ in general position \cite{pathflipsconv} and for generalized double circles \cite{pathflips}.

\section{Slide Graph of Caterpillars in Convex Position}

In this section, let $S$ be a set of $n$ points in convex position. Our goal is to prove that $G_{\mathcal{C}}^{\text{slide}}(S)$ is connected. We start with the following lemma.

\begin{lemma}\label{Lem: CattostarConvex}
    Let $C\in \mathcal{C}(S)$ be a caterpillar and let $s$ be one of the endpoints of its spine. Then $C$ can be transformed into $\starat(s)$ using at most $n-1-\deg(s)$ slides.
\end{lemma}

\begin{figure}
    \centering
    \includegraphics[scale=0.6]{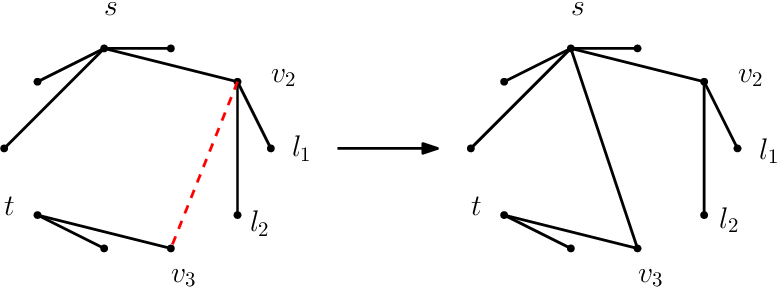}
    \caption{Edge $v_2v_3$ blocks leaves $l_1,l_2$ from being slid to $s$. Sliding the edge $v_2v_3$ results in caterpillar with spine $v_2sv_3t$, which allows $l_1,l_2$ to be slid to $s$.}
    \label{Fig: edgeblock}
\end{figure}

\begin{proof}
    Recall that $s$ is not a leaf vertex in $C$, and therefore has degree at least $2$. Assume that the spine of $C$ is the path $s,v_2,\dots,v_{k}$. For each $i\in \{2,\dots, k\}$, starting with $v_2$, first slide all of the leaves attached to $v_i$ visible from $s$ to $s$, starting from the one closest to $s$. Then, slide edge $v_iv_{i+1}$ along the edge $sv_i$. Since all previously slid edges form a star, $s$ is an endpoint of the spine, and the rest is part of the original caterpillar, this is a valid slide (see Figure \ref{Fig: edgeblock}). Now, at least one of the leaves attached to $v_i$ is visible from $s$, so slide it to $s$. Repeat this until $v_i$ becomes a leaf. In the entire procedure, each edge is slid exactly once except for the edges at $s$ in the initial caterpillar $C$, thus we need at most $n-1-\deg(s)$ slides. $\qed$
\end{proof}

We now observe that for all $u,v \in S$, $\starat(u)$ and $\starat(v)$ can be transformed into one another using $n-2$ slides. Then, Theorem \ref{Thm: Convexpillars} follows from this fact together with Lemma \ref{Lem: CattostarConvex}. Given two caterpillars $C_1$, $C_2$ on $S$, we first transform $C_1$ to a star in at most $n-3$ steps, then one star into another in at most $n-2$ steps and lastly we transform the star into $C_2$ in at most $n-3$ more steps. We get the following result as a corollary.

\begin{corollary}\label{Cor: slidingpaths}
    Let $n\ge 5$ and $P,Q$ be two plane spanning paths on $S$. Then, $P$ can be transformed into $Q$ using at most $2n-6$ slides in $G_{\mathcal{C}}^{\text{slide}} (S)$.
\end{corollary}

\begin{proof}
       Since $n\ge 5$, there exists a vertex $v$ of $S$ which has degree at least $2$ in both $P$ and $Q$. Then, if we split both $P$ and $Q$ at $v$ and perform the algorithm from Theorem \ref{Thm: Convexpillars} on all $4$ subpaths, we will transform both $P$ and $Q$ to $\starat(v)$ in at most $n-3$ moves each. Thus, since all moves are reversible we can slide from $P$ to $\starat(v)$ and finally to $Q$ in $2n-6$ moves. We can choose $v$ arbitrarily since in a path there are no leaves that are not head/tail-leaves. $\qed$
\end{proof}

The upper bound on the diameter of $G_{\mathcal{C}}^{\text{slide}}(S)$  that we obtain in Theorem \ref{Thm: Convexpillars} is most likely not tight. As mentioned before, it is known that for slides in plane trees, there is a lower bound of $1.5n-5$ \cite{slidelb}. Even further, the trees that achieve this lower bound are caterpillars, thus the gap between the lower and upper bound is big. Thus, it would be interesting to find better lower bounds for all of the reconfiguration graphs for $S$ in convex position.

\section{Reconfiguration Graphs of Caterpillars in General Position}

\begin{figure}
    \centering
    \includegraphics[scale=0.5]{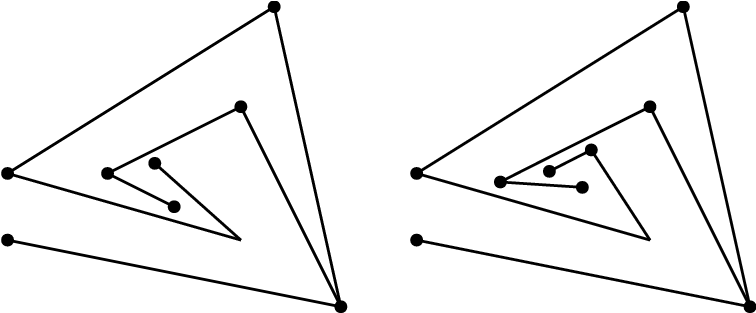}
    \caption{Caterpillars on $9$ and $10$ points with no available slides.}
    \label{Fig: noflips}
\end{figure}

In this section, let $S$ be a set of $n$ points in general position in the plane.  
As mentioned in the introduction, we can show that if $n\ge 9$, then $G_{\mathcal{C}}^{\text{slide}}(S)$ has at least one isolated vertex. We give examples for the cases of $n=9,10$ in Figure \ref{Fig: noflips}. The construction of the examples can easily be generalized to larger point sets.

Now, we move on to the proof of Lemma \ref{Thm: starsanddoublestars}. We will instead prove the following stronger statement which implies Lemma \ref{Thm: starsanddoublestars}.

\textbf{Statement A:} \textit{ Assume that $C$ is a plane spanning double-star on a point set $S$ such that $|S|=n$. Let $u$, $v$ be the centers of $C$. Then there exists a slide sequence in $G_{\mathcal{C}}^{\text{slide}}(S)$ from $C$ to $\starat(v)$. Moreover, $u$ is a spine vertex or a head/tail-leaf at every step of the sequence.}
 
\begin{figure}
    \centering
    \includegraphics[scale=0.55]{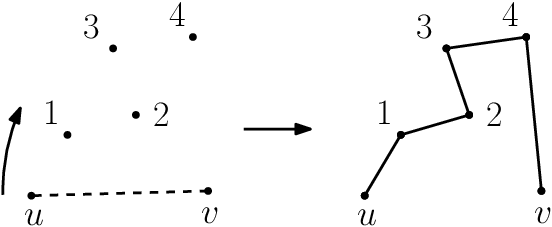}
    \caption{Radial $uv$ path.}
    \label{Fig: Radial}
\end{figure}

For a set of points $S$ and two adjacent vertices $u,v$ on convex hull of $S$, we define the \emph{radial} $uv$ path as the spanning path between $u$ and $v$ which visits the vertices of $S$ in the order they are seen when rotating the line through $uv$ around $v$ towards the interior of $S$, see Figure \ref{Fig: Radial}. 

\textbf{Statement B:} \textit{Let $u,v$ be two adjacent vertices of the convex hull of a point set $T$. Let $P$ be a radial $uv$ path in $T$. Then, there is a slide sequence in $G_{\mathcal{C}}(T)$ which transforms $P$ to a caterpillar whose spine uses all of the edges of the convex hull of $T$ except the edge $uv$. Moreover, $u$ is a spine vertex or a head/tail-leaf at every step of the sequence.}

To prove Statement A for a point set $S$ of $n$ points, we will assume that Statement B is true for any point set $T$ of size at most $n-1$. Then we will prove correctness of Statement B, using Statement A for point sets of size at most $n-2$. 
We now prove Statement A.

\begin{figure}[h]
    \centering
    \includegraphics[scale=0.5]{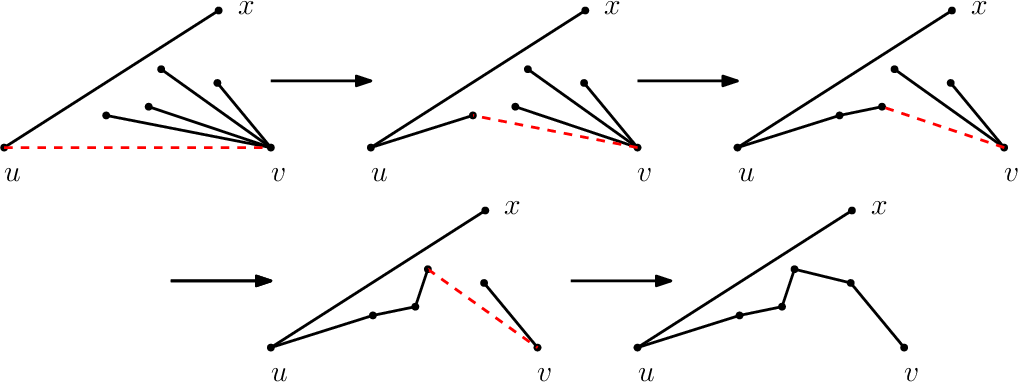}
    \caption{Constructing the radial path from $u$ to $v$.}
    \label{Fig: makingradialpath}
\end{figure}

\begin{proof}[of Statement A]

If $n=2$, there is nothing to do, so we assume that we can construct a slide sequence for all point sets of size at most $n-1$. As long as we can slide an edge connecting $u$ to a leaf so that it is attached to $v$, we do it. After some number of steps, we reach a leaf $x$ such that the triangle spanned by $u$, $x$ and $v$ contains some leaves $l_1,l_2\dots l_j$ (ordered as seen from $v$ when rotating $uv$ towards $x$) attached to $v$ in the interior of the triangle, as in the first part of Figure \ref{Fig: makingradialpath}. But now, we can apply Statement A inductively on the set of points inside the triangle spanned by $u$, $x$, and $v$ and make a star at $v$ on this point set without altering any of the edges outside of it. The only problematic case is if $u$, $x$, and $v$ span the entire convex hull of $S$.  In this case, we use the sequence of slides which can be seen in Figure \ref{Fig: makingradialpath}, which constructs a radial path $P$ from $u$ to $v$. Let $C_0$ be the resulting caterpillar, whose spine contains $P$ as a subpath. 

We now apply Statement B on path $P$ and point set $T$ spanned by $P$ to obtain a caterpillar $C \in \mathcal{C}$. Note that the sequence of slides transforming $P$ to $C$ in $G_{\mathcal{C}}^{\text{slide}}(T)$ kept $u,v$ as either endpoints of the spine or head/tail-leaves at every step, and hence it is a valid sequence of slides transforming $C_0$ to a caterpillar $C_1$ with $C$ as a subcaterpillar in $G_\mathcal{C}^{\text{slide}}(S)$.
 
Now, we can slide the edge $ux$ along the spine of $C$ to $v$. Then, we reverse all of the steps that created $C$ and $P$ and thus prove the theorem. It is important to note that we never create a tree that is not a caterpillar since our induction hypothesis preserves $u$ and $v$ as spine vertices. Even when we create the path $P$, we add extra edges of the spine between $u$ and $v$, but we still maintain the property that $u$ and $v$ are vertices of the spine. $\qed$
\end{proof}

Now we prove Statement B.

\begin{figure}[h!]
    \centering
    \includegraphics[scale=0.5]{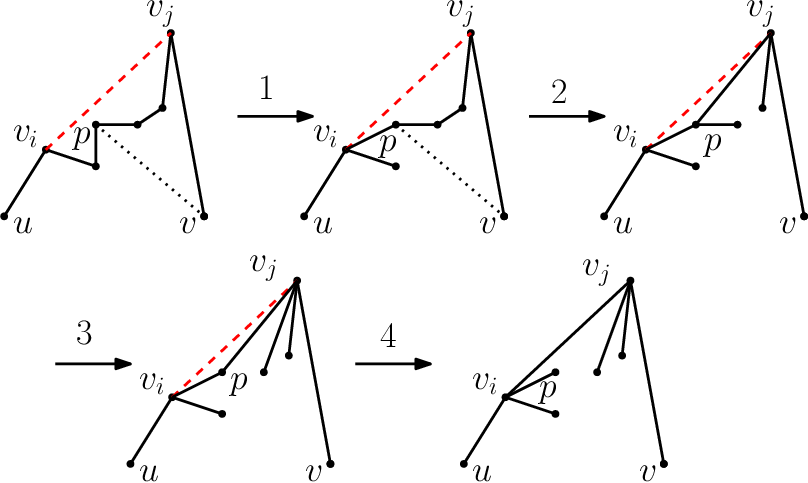}
    \caption{The process of adding the edge $v_iv_j$ to the caterpillar in proof of Statement B}
    \label{Fig: triplestar}
\end{figure}

\begin{figure}[h!]
    \centering
    \includegraphics[scale=0.6]{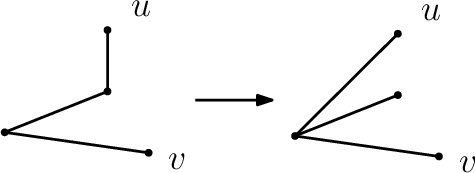}
    \caption{Base case for proof of Statement B}
 \label{Fig: basecase}
\end{figure}

\begin{proof} [of Statement B]
We prove Statement B by induction on $T$ with base case $T=4$ seen in Figure \ref{Fig: basecase}. We write $P=u,v_1,v_2,\dots, v_k,v$. Let $i,j$ be the smallest numbers such that $v_iv_j$ is an edge of the convex hull of $T$ which is not in $P$. We consider the point set $S_i^j= \{v_i, v_{i+1}\dots, v_j,v\}$. Let $p$ be the point inside the convex hull of $S_i^j$ which is visible from both $v_i$ and $v_j$ and closest to the segment $v_iv_j$. Then the path $v_i,\dots,p$ is a radial path from $v_i$ to $p$. By applying the reverse of the procedure in Figure \ref{Fig: makingradialpath} to the path $v_i,\dots,p$, it can be transformed into a double star with centers $v_i$ and $p$, see Step $1$ in Figure \ref{Fig: triplestar}. Again $p,\dots,v_j$ is a radial path from $p$ to $v_j$, so we transform $p,\dots,v_j$ into a double star with centers $p,v_j$ in exactly the same way, see Step $2$ in Figure \ref{Fig: triplestar}.  Now, by inductively applying Statement A on the set of points inside the triangle $v_i,p,v$, we can transform the double star with centers $v_i,p$  into a star with center $v_i$. Then, we apply Statement A to transform the double star at $p,v_j$ to the star at $v_j$, we can do this since $v_j,p$ are consecutive along the convex hull of $S_i^j\setminus \{v_i\}$. Lastly, we slide the edge $v_j p$ to $v_i$. See Steps 3 and 4 in Figure \ref{Fig: triplestar}, respectively. We repeat this entire procedure for every edge of the convex hull of $S'$ which is not in $P$. $\qed$
\end{proof}

\noindent With some extra work, we can drop the requirement that $u$ and $v$ are on the convex hull of $S$. 

\begin{corollary}\label{Cor: stars}
    Let $u$ and $v$ be arbitrary vertices of $S$. Then both $\starat(u)$ and $\starat(v)$ are connected to any double star with centers at $u$ and $v$.
\end{corollary}
\begin{proof}
    The line through $u$ and $v$ divides $S$ into two point sets $T,T'$ such that $u$ and $v$ are adjacent on the convex hull of both $T$ and $T'$. We then apply Lemma \ref{Thm: starsanddoublestars} on the double star with centers $u,v$ in $T$ and $T'$. $\qed$
\end{proof}

 Using Corollary \ref{Cor: stars}, we prove Theorem \ref{Thm: wellseparated}.

\begin{proof}[of Theorem \ref{Thm: wellseparated}]
    We first prove that $C$ can be transformed into a star. Denote the vertices of the spine of $C$ by $v_1,v_2,\dots, v_k$. We proceed by induction on $k$. The case $k=1$ is covered by Corollary \ref{Cor: stars}. Now, as $C$ is well-separated, we know that the convex hull of $S_{1,2}$ is disjoint from the rest of $S$. Thus we can apply Corollary \ref{Cor: stars} to $S_{1,2}$ and reduce the length of the spine of $C$ by one. Then the result follows by induction and Corollary \ref{Cor: stars}. $\qed$
\end{proof}

Before proving Theorem \ref{thm: sliidesconnnected}, we need the following slight strengthening of Corollary \ref{Cor: stars}.

\begin{lemma}
    Let $S$ be a set of points in the plane and $C$ a plane spanning caterpillar on $S$ with spine consisting of three vertices. Then, there is a slide sequence transforming $C$ to a star. 
\end{lemma}\label{Lem: triplestartostar}

\begin{figure}
    \centering
    \includegraphics[scale=0.5]{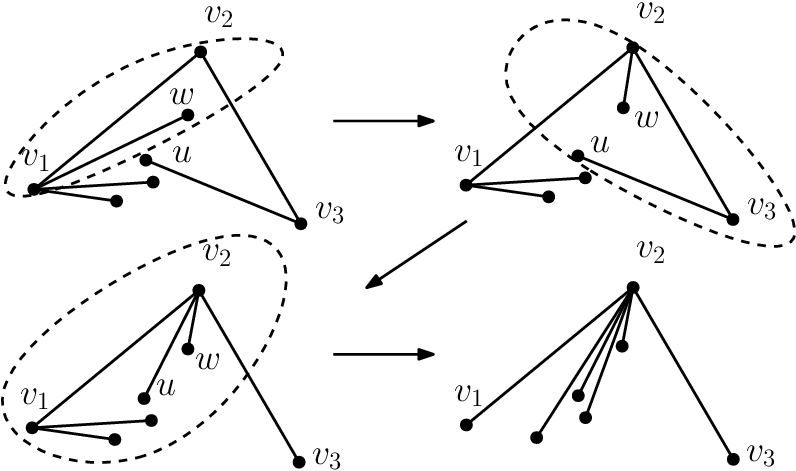}
    \caption{Visualization of the proof of Lemma \ref{Lem: triplestartostar}. At every step, we use Lemma \ref{Thm: starsanddoublestars} on the point set inside the dashed line.}
    \label{fig:3starsto2stars}
\end{figure}

\begin{proof}
    Let $v_1$, $v_2$ and $v_3$ be the three spine vertices of $C$, where $v_2$ is the central one. We will show that there exists a sequence of slides which shortens the spine to the vertices $v_2$ and $v_3$. Then, by Corollary \ref{Cor: stars}, $C$ will be connected to both $\starat(v_2)$ and $\starat(v_3)$. Consider the set $S_{1,2}$. If its convex hull contains no points connected to $v_3$, we are done. Otherwise, let $u$ be the first such point we encounter when rotating the edge $v_2v_3$ with $v_3$ fixed. If the triangle defined by $\{u,v_2,v_3\}$ contains no points connected to $v_1$, by Lemma \ref{Thm: starsanddoublestars}, we may slide $uv_3$ to $uv_2$. Otherwise, let $w$ be such a point. Since $u$ was chosen to be minimal, the triangle of $\{w,v_1,v_2\}$ contains no points connected to $v_3$, and thus by Lemma \ref{Thm: starsanddoublestars} we can slide $wv_1$ to $wv_2$, see Figure \ref{fig:3starsto2stars} for the visualization of the entire process. By repeating this process, we are left with no points connected to $v_3$ inside the convex hull of $S_{1,2}$ and thus can slide all leaves connected to $v_1$ into leaves connected to $v_2$, finishing the proof. $\qed$
    
\end{proof}

\begin{proof}[of Theorem \ref{thm: sliidesconnnected}]
    If $n\ge  9$, $G_{\mathcal{C}}^{\text{slide}}(S)$ is disconnected by examples seen in Figure \ref{Fig: noflips}. If $n\le 6$ then every tree is a caterpillar so the result follows from \cite{slideconn}. If $n=7$ we checked the results computationally \cite{github}. To make our computations more efficient we used Lemma \ref{Lem: triplestartostar} which tells us that we only need to check caterpillars with spine of at least $4$ vertices and that we could stop computations as soon as we decreased the length of spine by one. If $n=8$ we used the same approach but we found an isolated vertex. $\qed$
\end{proof}

To end the section, we now focus on rotations and prove Proposition \ref{prop: rotations}.

\begin{figure}[h]
    \centering
    \includegraphics[scale = 0.6]{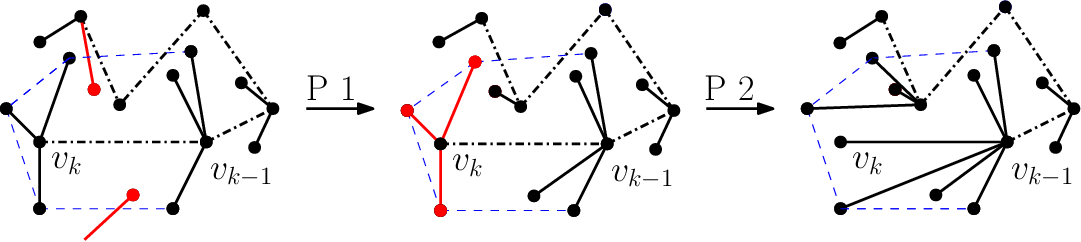}
    \caption{Algorithm described in the proof of Proposition \ref{prop: rotations}, spine edges are dash-dotted in the picture. Each phase lasts until all of the red edges are rotated.}
    \label{Fig:prop1}
\end{figure}

\begin{proof}[of Proposition \ref{prop: rotations}]
    We will prove that for any caterpillar $C\in \mathcal{C}(S)$ we can find a sequence of rotations transforming $C$ into a star. We will prove this by induction on $k$, where $k$ is the length of the spine of $C$. Case $k=2$ follows from Corollary \ref{Cor: stars} as every slide is a rotation. Now assume that the statement holds for $k-1$ and let $C$ be a caterpillar with spine $v_1,v_2,\dots, v_k$. If the convex hull of $S_{k-1,k}$ is disjoint from the rest of $S$, the result follows by inductive hypothesis and Corollary \ref{Cor: stars}. Thus, we may assume that there are some vertices of $S\setminus S_{k-1,k}$ inside the convex hull of $S_{k-1,k}$. Let $I$ be the set of all such vertices. We now divide the algorithm into two phases, see Figure \ref{Fig:prop1}. 

    \textbf{Phase 1:} Let $L\subseteq I$ be the subset of $I$ such that a vertex is in $L$ if and only if it is a leaf connected to a vertex of the spine which is outside the convex hull of $S_{k,k-1}$. We repeat the following procedure as long as $L\neq \emptyset$. Find a leaf vertex $l$ in $L$ which sees at least one spine vertex $v_j$ lying inside the convex hull $S_{k,k-1}$ ($v_k$ and $v_{k-1}$ included). Then make a rotation so that $l$ is connected to $v_j$.

    \textbf{Phase 2:} Repeat the following procedure until $v_k$ has no leaves attached to it and thus becomes a head-leaf itself. Pick a leaf $l$ attached to $v_k$ such that $l$ sees at least one spine vertex $v_j$ in $I$, such that $j\neq k$. After Phase 1 is completed such a leaf always exists. Now make a rotation so that $l$ is connected to $v_j$.

    See Figure \ref{Fig:prop1} for a visualization of both phases.
    
    After completing both phases the spine of $C$ is shortened by one and the result follows by induction. $\qed$
\end{proof}

\section{Connected Components in the Flip Graph of Plane Paths}\label{Sec: flips}

\begin{figure}[h]
    \centering
    \includegraphics[scale = 0.6]{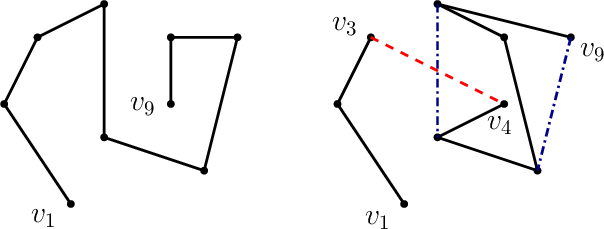}
    \caption{The path on the left is a generalized peeling path, the path on the right is not due to the edge $v_3v_4$ crossing into the convex hull of $\{v_4,\dots, v_9\}$.}
    \label{Fig: Generalizedpeels}
\end{figure}

In this section, we shift our focus and consider the flip graph of plane spanning paths, particularly we study the sizes of connected components of this graph.  Recall that, for a point set $S$ in general position, we denote the set of all plane spanning paths by $\mathcal{P}(S)$ and the corresponding flip graph by $G_{\mathcal{P}}(S)$. We say that a plane spanning path $P=v_1,\dots, v_n$ is a \emph{generalized peeling path} if $v_1$ is a vertex of the convex hull of $S$ and for $i>1$, the vertex $v_i$ is a vertex of the convex hull of $S_i = S\setminus \{v_1,\dots, v_{i-1}\}$ and it is visible from $v_{i-1}$ (so, no edge of the convex hull of $S_i$ blocks it). For an illustration, see Figure \ref{Fig: Generalizedpeels}. Obviously, this definition is dependent on the way we orient the path. For example, if we reverse the path on the left of Figure \ref{Fig: Generalizedpeels} it is not a generalized peeling path. It is not hard to see that for each point set $S$, there are at least $2^{n-1}$ generalized peeling paths. This is because at every step of constructing a generalized peeling path, we have at least two choices for the next vertex along the path (except at the last step). Therefore, to prove Theorem \ref{Thm: largecomponent} we only need to prove that generalized peeling paths lie in a single component of $G_\mathcal{P}(S)$.

As mentioned before, we move the proof of Theorem \ref{Thm: largecomponent} to the appendix. Instead, we will use the generalized peeling paths to estimate the number of well-separated caterpillars.

One can easily observe that any generalized peeling path on $k$ vertices is also a well-separated caterpillar on $k$ vertices when the orientation is reversed. Further, if $S$ is a point set on $n$ points and $S' \subseteq S$ is a subset of $k\le n$ points then it is easily seen that any generalized peeling path on $S'$ defines at least one well-separated caterpillar on $S$, again with spine having reverse orientation of the original path. Therefore, we get a lower bound on the number of well-separated caterpillars as follows. If $\text{GP}(k)$ is the number of generalized peeling paths on $k$ vertices and $\text{WS}(n)$ is the number of well-separated caterpillars on $n$ vertices, then

$$ \text{WS}(n) \ge \sum_{k=1}^n \binom{n}{k}\text{GP}(k)\ge  \sum_{k=1}^n \binom{n}{k}2^{k-1} = \frac{1}{2}(3^n-1).$$
Therefore, we have proven that $G_{\mathcal{C}}^{\text{slide}}(S)$ has a connected component of size at least $\frac{1}{2}(3^n-1)$.

 In order to prove Theorem \ref{thm: smallcomponents} we first need to collect some results regarding the flip graph of all plane spanning paths with a fixed endpoint $u\in S$, denoted by $G_{\mathcal{P}}^u(S)$. 

\begin{lemma} \label{Lem: girth}
    The graph $G_{\mathcal{P}}^u(S)$ has girth of at least six.
\end{lemma}

\begin{proof}
         In $G_{\mathcal{P}}^u(S)$ any flip that changes both endpoints of a path is forbidden. Thus, any allowed flip can be viewed as a suffix reversal of the path, changing a path of the form $u,v_2,\dots, v_n$ to $u,v_2,\dots,v_{k-1},v_n,v_{n-1},\dots, v_k$ for a choice of $v_k$ which keeps the resulting path plane. Therefore, $G_{\mathcal{P}}^u(S)$ can be considered as a subgraph of the corresponding suffix reversal graph (also known as the pancake graph), which is well known to have girth 6~\cite{GirthOfPancakeGraphs}. $\qed$
\end{proof}

\begin{lemma} \label{Lem: no isolated paths}
    The graph $G_{\mathcal{P}}^u(S)$ has no isolated vertices.
\end{lemma}

\begin{proof}
        Assume that $P=u,v_2,\dots, v_{n-1},v_n$ is isolated in $G_{\mathcal{P}}^u(S)$. This means that $v_n$ sees no vertex of $S$ other than $v_{n-1}$. However, the first vertex seen when rotating $v_nv_{n-1}$ towards $v_{n-2}$ is either on the interior of the triangle $v_n,v_{n-1},v_{n-2}$ or $v_{n-2}$, and is therefore always visible from $v_n$. $\qed$

\end{proof}

The following lemma characterizes vertices of degree one in $G_{\mathcal{P}}^u(S)$. The proof is technical and involved so we move it to the appendix. 
\begin{lemma}\label{Lem: Paths of degree one}
    If $P=u,v_2,\dots,v_n$ is a path of degree one in $G_{\mathcal{P}}^u(S)$, then $v_{n-1}$, $v_n$ and $v_{n-2}$ are consecutive vertices of the convex hull of $S$ and the interior of the triangle $v_{n-1},v_n,v_{n-2}$ is disjoint from $S$.
 \end{lemma}

 \begin{figure}[ht]
     \centering
     \includegraphics[scale =0.5]{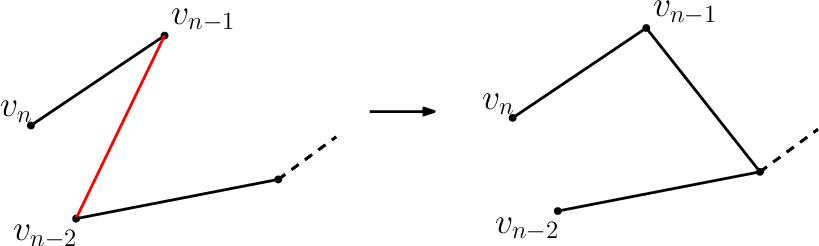}
     \caption{Path of degree one and effect of the single flip on it.}
     \label{Fig: paths of degree one}
 \end{figure}

We put everything together in the following proposition. 

\begin{proposition}\label{prop: smallcomps}
    Let $S$ be a set of points in general position and $u$ a vertex of $S$. Then $G_\mathcal{P}^{u}(S)$ does not have a connected component on at most $4$ vertices.
\end{proposition}

\begin{proof}
       By Lemma \ref{Lem: no isolated paths}, the graph $G_{\mathcal{P}}^u(S)$ has no isolated vertices. Further, if $P=u,\dots,p_{n-1},p_n$ and $Q= u,\dots, q_{n-1},q_n$ are two distinct paths of degree one in $G_{\mathcal{P}}^u(S)$ then they admit a single flip which transforms them to $P'=u,\dots,p_{n},p_{n-1}$ and $Q'=u,\dots,q_n,q_{n-1}$ respectively. And as $P\neq Q$, it follows that $P'\neq Q'$.   This forbids $K_{1,3},P_2$ and $P_3$ as connected components of $G_{P}^u(S)$, see Figure \ref{Fig: paths of degree one}. By Lemma \ref{Lem: girth}, cycles of length $3$ and $4$ are also forbidden. Therefore, the only remaining graph on $4$ vertices that may be a connected component is $P_4$. The proof of this fact is rather technical and moved to the appendix. $\qed$
\end{proof}

Finally, we prove Theorem \ref{thm: smallcomponents}. 

\begin{proof}[of Theorem \ref{thm: smallcomponents}]
    Assume that $G_{\mathcal{P}}(S)$ contains a connected component on $7$ vertices. Call these vertices $P_1,P_2,\dots P_7$. Consider $P_1$, a plane spanning path on $S$ with endpoints $u,v$. By Proposition \ref{prop: smallcomps}, we know that if we fix $u$, $P_1$ is still connected to $4$ vertices via flips in $G_{\mathcal{P}}^u(S)$. If any of these $4$ is not one of $P_2,P_3,\dots P_7$ we get a component of size at least $8$, contradicting our assumption. Therefore we may assume that those $4$ vertices are $P_2,P_3,P_4,P_5$, which all need to have $u$ as one of their endpoints. By a symmetric argument, we see that $4$ more of the vertices $P_2,P_3,\dots P_7$ need to have $v$ as one of their endpoints. Thus by the pigeonhole principle, we can assume that $P_1,P_2,P_3$ all have $u,v$ as endpoints.  But now we can apply the identical argument to $P_4$ which has endpoints $u,w$ where $w\neq v$. Therefore in the $7$ paths $P_1,P_2,\dots, P_7$ at least $3$ have $u,v$ as endpoints, at least $3$ have $u,w$ as endpoints and for each  $x\in \{u,v,w\}$ at least $4$ of the paths have $x$ as one the endpoints. This is clearly impossible so the result follows. $\qed$
\end{proof}

As a possible step towards resolving Conjecture \ref{conj: main}, it would be interesting to determine if similar observations about the structure of $G_{\mathcal{P}}(S)$ and $G_{\mathcal{P}}^u(S)$ can be used to show nonexistence of larger sized connected components in $G_{\mathcal{P}}(S)$.

\section*{Acknowledgements}
We would like to thank Jan Kyn\v{c}l and Pavel Valtr for proposing this problem to us in the Combinatorial Problems Seminar at Charles University where this work started and for all of their advice and help during the work on this problem. We also thank Daniel Perz for pointing out the existence of an isolated vertex in $G_\mathcal{C}^{\text{slide}}(S)$, for his help in finding literature for this paper and helpful comments.



\bibliographystyle{plain}
\bibliography{bibliography.bib}

\appendix

\section{Proof of Theorem \ref{Thm: largecomponent}}
     First, we show that for any vertex $u$ of the convex hull $S$, all generalized peeling paths starting at $u$ are in the same component of $G_{\mathcal{P}}^u(S)$. Clearly, the result holds for $n=3$ so we assume that $n>3$ and proceed by induction. Let $P=u,p_2, \dots, p_n$ and $Q=u,q_2, \dots, q_n$ be generalized peeling paths. If $p_2=q_2$, the result follows by induction so we assume $p_2 \neq q_2$. We construct a generalized peeling path $F$ on $S\setminus \{u\}$ with endpoints $p_2,q_2$ as follows. Let $C$ be the convex hull of $S \setminus (\{u\} \cup \mathcal{F})$ where $\mathcal{F}$ is the set of points currently added to $F$. Let $p$ be the last point added to $F$. Choose the longer of the two paths from $p$ to $q_2$ along the boundary of $C$ and add those vertices to $F$. Repeat this process until $C=\{q_2\}$ and then  add $q_2$ to $F$. Now, by induction, we can transform $P\setminus\{u\}$ to $F$, and then by a single flip we can transform the path $u,F$ to a path with starting edge $uq_2$. Again, by induction, we can transform this new path to $Q$. This proves that all generalized peeling paths starting at a vertex $u$ are in the same connected component of $G_{\mathcal{P}}^u(S)$. Further, for any other vertex $u'$ on the convex hull of $S$, we can construct a generalized peeling path from $u$ to $u'$ as before, thus finishing the proof. 

\section{Proof of Lemma \ref{Lem: Paths of degree one}}
    
\begin{figure}[ht]
     \centering
     \includegraphics[scale = 0.5]{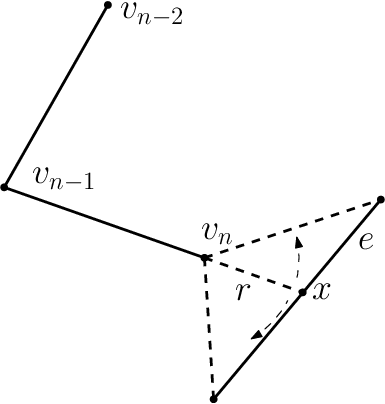}
     \caption{Illustration of the last part of the proof of Lemma \ref{Lem: Paths of degree one}. Since all of the dashed segments are uncrossed by edges of $P$, the first vertex seen in the triangles bounded by dashed segments plus $e$ is always visible from $v$.}
     \label{Fig: paths of degree one proof}
 \end{figure}

     Assume that $v_{n-2},v_n,v_{n-1}$ are not consecutive vertices of the convex hull. First note that $v_n$ sees at least one vertex in the triangle $v_n,v_{n-1},v_{n-2}$, as in the proof of Lemma \ref{Lem: no isolated paths}. Assume that it sees only one vertex inside the triangle as otherwise we are done. Therefore we want to prove that it sees a vertex outside the triangle. Let $e$ be the edge of $P$ whose intersection $x$ with the ray $r$ from $v_{n-1}$ through $v_n$ is closest to $v_n$. If $e$ is well defined (if there is at least one edge which crosses $r$), then $v_n$ sees the first vertex inside the triangle formed by $v_n,x$ together with the endpoints of $e$ which is encountered when rotating the segment $v_nx$ around $v_n$ in either direction, see Figure \ref{Fig: paths of degree one proof}. Otherwise, $v_n$ sees the first vertex encountered when rotating $r$ around $v_n$ towards $v_{n-2}$. If this vertex is $v_{n-2}$, then we are in one of two situations. Either all of the vertices of $S$ are inside the triangle $v_{n-1},v_n,v_{n-2}$, in which case $v_n$ must see at least two of them. Otherwise, all of the points of $S$ are on the opposite side of $v_nv_{n-1}$ from $v_{n-2}$. In this case, $v_n$ sees the first vertex encountered when rotating $r$ around $v_n$ away from $v_{n-2}$.

\section{Proof that $P_4$ cannot be a connected component of $G_{\mathcal{P}}^u(S)$}

\begin{figure}
    \centering
    \includegraphics[scale =0.6]{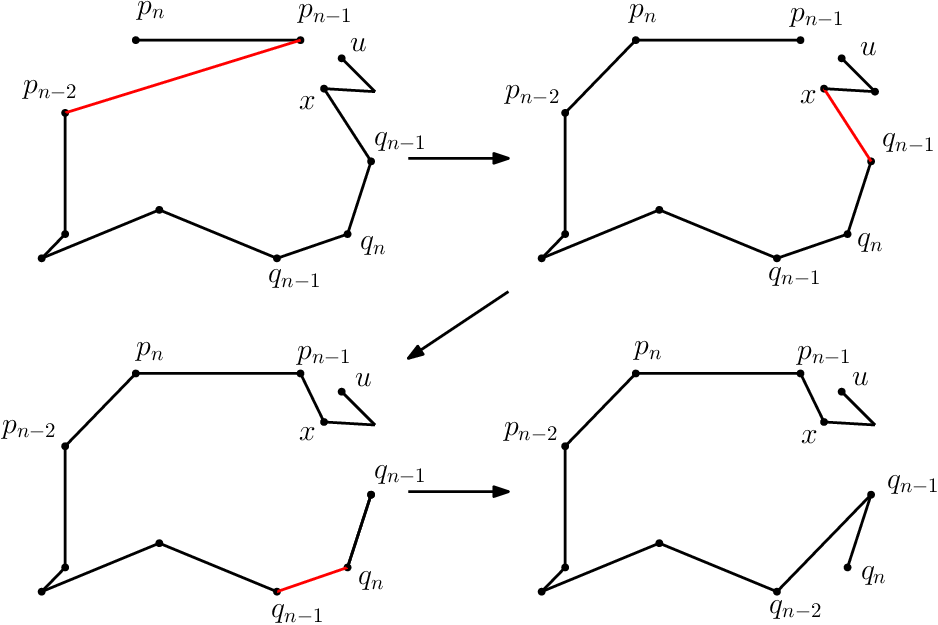}
    \caption{Transforming a path of degree one into another path of degree one using three flips. Intermediate paths have degree more than two and hence this is not the entire connected component.}
    \label{Fig: pathoflengthfour}
\end{figure}

Recall that $|S|\ge 5$. Let $P=u,p_2,\dots, p_n, Q= u,q_2,\dots,q_n\in G_{\mathcal{P}}^u(S)$ be two paths of degree one such that $p_{n-1},p_n,p_{n-2}$ and $q_{n-1},q_n,q_{n-2}$ are (not necessarily disjoint) triples of consecutive vertices of the convex hull of $S$. If $P$ can be transformed into $Q$ in three flips, such that none of the intermediate paths have degree larger than two, then we would have a connected component $P_4$. Intermediate paths are $P' = u,\dots, p_{n-2},p_n,p_{n-1}$ and $Q' =u,\dots, q_{n-2},q_n,q_{n-1}$. Assume that $p_{n}\neq q_{n}$. Let $x$ be the last vertex between $u$ and $q_{n-1}$ in path $P'$. If $P'$ and $Q'$ are connected by a single flip, segments $xq_{n-1}$ and $xp_{n-1}$ are uncrossed by an edge of $Q'$ nor $P'$. Therefore, $p_{n-1}$ sees vertices $p_{n-2}$ and $x$. If $x$
is not adjacent to $p_{n-1}$ along the convex hull, then we are done as $p_{n-1}$ always sees the vertices adjacent to itself on the convex hull, see Figure \ref{Fig: pathoflengthfour} . If $x$ is adjacent to $p_{n-1}$ on the convex hull. Let $y$ be the edge preceding $x$ in $P'$. We know that $p_{n-2}\neq y \neq x$ as otherwise $x$ would lie between $p_{n-1}$ and $q_{n-1}$ in $P'$. Thus, $p_{n-1}$ will see the first vertex encountered when rotating segment $p_{n-1}x$ towards $y$. Intermediate paths have degree more than two and hence this is not the entire connected component.

If $p_{n}=q_{n}$, the situation is slightly different. In this case, $p_{n-2}=q_{n-1}$ and $p_{n-1}=q_{n-2}$ so we need to consider the line segments $p_{n-1}p_{n-3}$ and $p_{n-2}p_{n-3}$. If the triangle $p_{n-1},p_{n-2},p_{n-3}$ is empty, then rotating the segment $p_{n-2}p_{n-3}$ around $p_{n-2}$ inside the path $Q'=u,p_1,\dots,p_{n-1},p_n,p_{n-2}$ allows $p_{n-2}$ to see three vertices: $p_{n-1}$, $p_{n-3}$ and the first vertex encountered in the rotation. Therefore $P$ and $Q$ cannot be vertices of a connected component isomorphic to $P_4$. 

 \end{document}